\documentclass[12pt,a4paper]{article}
\usepackage{amsmath}
\usepackage{amsthm}
\usepackage{amsfonts}
\usepackage{amssymb}
\usepackage{latexsym}

\usepackage[usenames]{color}
\definecolor{light-blue}{rgb}{0.8,0.85,1}
\definecolor{light-gray}{rgb}{0.8,0.8,0.5}
\definecolor{light-red}{rgb}{1,0.7,0.7}

\newtheorem{remark}{Remark}
\newtheorem{construction}{Construction}

\theoremstyle{definition}
\newtheorem{definition}{Definition}
\newtheorem{example}{Example} % The '*' makes it unnumbered

\begin{document}

\title{\Large{\textbf{Prolongations of quasigroups}}}
\author{\normalsize {V.A. Shcherbacov}}

 \maketitle

\begin{abstract}
\noindent
We give an overview of various prolongations of quasigroups. Two step prolongation procedure is proposed.

\medskip

\noindent \textbf{2000 Mathematics Subject Classification:} 20N05, 05B15

\medskip

\noindent \textbf{Key words and phrases:} quasigroup, quasigroup prolongation
\end{abstract}

\bigskip
\bigskip\bigskip

\hfill \begin{parbox}{100mm}{\textbf{In memoriam: Galina Borisovna  Belyavskaya (1940 - 2015).}}
\end{parbox}
\bigskip
\bigskip

Given information was prepared for submission  at the conference LOOPS'15 (28 June - 04 July, 2015,  Ohrid, Macedonia).

\section{Quasigroup prolongations} \label{QUAS_PROLONG}

Quasigroup prolongation is quit natural way of construction of a finite quasigroup of order $n+k$ ($k\leq n$) from a quasigroup of order $n$.
We start from some definitions. Mainly we follow \cite{DK1, VD,BEL_EXTENSIONS, DER_DUD}.

R.H. Bruck proposed to use transversals ($n$-transversals) for prolongation of a quasigroup   \cite{6}.
We give Belousov construction of quasigroup prolongation. If transversal elements are situated on main diagonal, we obtain Bruck construction.

\begin{construction} \label{CONSTRUCTION_1}
We prolong the Latin square and quasigroup $(Q, \ast)$ of order $3$ to  Latin square and  corresponding quasigroup $(Q^{\prime}, \star)$ of order $4$ in the following way.

We add additional column from the right  and additional row below, transpose in these new cell all marked (transversal)  elements in their fixed order and fill all remaining  empty cells by the symbol "4".
\[
\begin{array}{lllll}
\begin{array}{lll}
{1} & {2} & \colorbox{green}{3} \\
\colorbox{green}{2} &{3} & 1 \\
3 &	\colorbox{green}{1} & {2} \\
\end{array}
& \rightarrow &
\begin{array}{llll}
1 &	2 &	\Box &	\colorbox{green}{3}\\
\Box  &	3 &	1 &	\colorbox{green}{2}\\
3&\Box	&	2&	\colorbox{green}{1}\\
\colorbox{green}{2} &	\colorbox{green}{1}&	\colorbox{green}{3} &	\Box\\
\end{array}
& \rightarrow &
\begin{array}{llll}
1 &	2 &	\colorbox{yellow}{4} &	\colorbox{green}{3}\\
\colorbox{yellow}{4} &	3 &	1 &	\colorbox{green}{2}\\
3&	\colorbox{yellow}{4}&	2&	\colorbox{green}{1}\\
\colorbox{green}{2} &	\colorbox{green}{1}&	\colorbox{green}{3} &	\colorbox{yellow}{4}\\
\end{array}
\end{array}
\]
\end{construction}
\iffalse
&
\begin{array}{l|llll}
\star & 1 & 2 & 3 & {4}\\
\hline
1&1 &	2 &	\colorbox{yellow}{4} &	\colorbox{green}{3}\\
2&\colorbox{yellow}{4} &	3 &	1 &	\colorbox{green}{2}\\
3&3&	\colorbox{yellow}{4}&	2&	\colorbox{green}{1}\\
4&\colorbox{green}{2} &	\colorbox{green}{1}&	\colorbox{green}{3} &	\colorbox{yellow}{4}\\
\end{array}
\fi

It is easy to see that the initial Latin square from Construction \ref{CONSTRUCTION_1} has more than one transversal. Using  colored boxes we "isolated" in the initial Latin square three disjoint transversals.
\begin{example} \label{Latin_squere_Of_Ordr_3}
\[
\begin{array}{lll}
\colorbox{yellow}{1} & \colorbox{green}{2} & \colorbox{light-blue}{3} \\
\colorbox{light-blue}{2} &\colorbox{yellow}{3} & \colorbox{green}{1} \\
\colorbox{green}{3} &	\colorbox{light-blue}{1} & \colorbox{yellow}{2} \\
\end{array}
\]
\end{example}

There exists a possibility to generalize Construction \ref{CONSTRUCTION_1} and to make  a quasigroup prolongation using $2, 3, \dots, n$ \textbf{disjoint} transversals \cite{Damm_DISS, Yamamoto}. We demonstrate  generalization of  Construction \ref{CONSTRUCTION_1}
on the following examples.

\begin{example} \label{Example_261}
We prolong the Latin square  of order $3$ to  Latin square  of order $5$ in the following way. \index{prolongation!Yamamoto} %\index{prolongation!Damm}
\end{example}

Step 1.
We add two additional columns from the right   and two additional rows below, transpose in these new cell all marked (transversal)  elements in their fixed order (i.e, we project these elements along rows and along columns).
Notice we can add these two columns and two rows in any suitable  place of a given Latin square.

\[
\begin{array}{lllll}
\Box & \Box & \colorbox{light-blue}{3} & \colorbox{yellow}{1} &  \colorbox{green}{2}  \\
\colorbox{light-blue}{2} & \Box & \Box & \colorbox{yellow}{3} &  \colorbox{green}{1}   \\
\Box & \colorbox{light-blue}{1} & \Box & \colorbox{yellow}{2} &  \colorbox{green}{3}   \\
\colorbox{yellow}{1} & \colorbox{yellow}{3} & \colorbox{yellow}{2} &  \Box & \Box \\
\colorbox{green}{3} & \colorbox{green}{2} & \colorbox{green}{1}  &  \Box & \Box  \\
\end{array}
\]

Step 2.

Fill all   empty after transposition transversal cells by the symbols "4, 5". In transversal cells of a fixed transversal we put the same element. Here yellow transversal we fill by the element "4".

\[
\begin{array}{lllll}
4 & 5 & \colorbox{light-blue}{3} & \colorbox{yellow}{1} &  \colorbox{green}{2}  \\
\colorbox{light-blue}{2} & 4 & 5 & \colorbox{yellow}{3} &  \colorbox{green}{1}   \\
5 & \colorbox{light-blue}{1} & 4 & \colorbox{yellow}{2} &  \colorbox{green}{3}   \\
\colorbox{yellow}{1} & \colorbox{yellow}{3} & \colorbox{yellow}{2} &  \Box & \Box \\
\colorbox{green}{3} & \colorbox{green}{2} & \colorbox{green}{1}  &  \Box & \Box  \\
\end{array}
\]

Step 3.
In remaining right bottom empty  square  we put a quasigroup of order $2$ defined on the set $\{ 4, 5\}$.

\[
\begin{array}{lllll}
4 & 5 & \colorbox{light-blue}{3} & \colorbox{yellow}{1} &  \colorbox{green}{2}  \\
\colorbox{light-blue}{2} & 4 & 5 & \colorbox{yellow}{3} &  \colorbox{green}{1}   \\
5 & \colorbox{light-blue}{1} & 4 & \colorbox{yellow}{2} &  \colorbox{green}{3}   \\
\colorbox{yellow}{1} & \colorbox{yellow}{3} & \colorbox{yellow}{2} &  5 & 4 \\
\colorbox{green}{3} & \colorbox{green}{2} & \colorbox{green}{1}  &  4 & 5  \\
\end{array}
\]

Modification of Step 1. We can change the order of colored rows and/or  columns.
\[
\begin{array}{lllll}
4 & 5 & \colorbox{light-blue}{3} & \colorbox{yellow}{1} &  \colorbox{green}{2}  \\
\colorbox{light-blue}{2} & 4 & 5 & \colorbox{yellow}{3} &  \colorbox{green}{1}   \\
5 & \colorbox{light-blue}{1} & 4 & \colorbox{yellow}{2} &  \colorbox{green}{3}   \\
\colorbox{green}{3} & \colorbox{green}{2} & \colorbox{green}{1}  &  5 & 4  \\
\colorbox{yellow}{1} & \colorbox{yellow}{3} & \colorbox{yellow}{2} &  4 & 5 \\
\end{array}
\]

\begin{remark}
It is clear that by prolongation of Latin squares and quasigroups we can situate additional columns and rows not only from the right and in the bottom of initial Latin square, but in any other suitable place.
\end{remark}

\begin{example}\label{Example_262}
We prolong the Latin square  of order $3$ to  Latin square of order $6$ in the following way.
\end{example}
Step 1.

We add three additional columns from the right  and three additional rows below, transpose in these new cell all marked (transversal)  elements in their fixed order.

\[
\begin{array}{llllll}
{\Box} & \Box & \Box & \colorbox{yellow}{1} &  \colorbox{green}{2} & \colorbox{light-blue}{3} \\
\Box & {\Box} & \Box & \colorbox{yellow}{3} &  \colorbox{green}{1} & \colorbox{light-blue}{2} \\
\Box & \Box & {\Box} & \colorbox{yellow}{2} &  \colorbox{green}{3} & \colorbox{light-blue}{1} \\
\colorbox{yellow}{1} & \colorbox{yellow}{3} & \colorbox{yellow}{2} & \Box & \Box & \Box \\
\colorbox{green}{3} & \colorbox{green}{2} & \colorbox{green}{1}  &  \Box & \Box & \Box \\
\colorbox{light-blue}{2} &	\colorbox{light-blue}{1} & \colorbox{light-blue}{3} &  \Box & \Box & \Box \\
\end{array}
\]

Step 2.

Fill all remaining  empty after transposition transversal cells by the symbols "4, 5, 6" in their  "old transversal order", i.e. we put the symbol $4$ in all empty cells of light-blue transversal.

\[
\begin{array}{llllll}
\colorbox{yellow}{6} & \colorbox{green}{5} & 4 & \colorbox{yellow}{1} &  \colorbox{green}{2} & \colorbox{light-blue}{3} \\
4 & \colorbox{yellow}{6} & \colorbox{green}{5} & \colorbox{yellow}{3} &  \colorbox{green}{1} & \colorbox{light-blue}{2} \\
\colorbox{green}{5} & 4 & \colorbox{yellow}{6} & \colorbox{yellow}{2} &  \colorbox{green}{3} & \colorbox{light-blue}{1} \\
\colorbox{yellow}{1} & \colorbox{yellow}{3} & \colorbox{yellow}{2} & \Box & \Box & \Box \\
\colorbox{green}{3} & \colorbox{green}{2} & \colorbox{green}{1}  &  \Box & \Box & \Box \\
\colorbox{light-blue}{2} &	\colorbox{light-blue}{1} & \colorbox{light-blue}{3} &  \Box & \Box & \Box \\
\end{array}
\]

Step 3.

In remaining right bottom empty  square  we put any quasigroup of order $3$ defined on the set $\{ 4, 5, 6 \}$.

\[
\begin{array}{llllll}
6 & 5 & 4 & \colorbox{yellow}{1} &  \colorbox{green}{2} & \colorbox{light-blue}{3} \\
4 & 6 & 5 & \colorbox{yellow}{3} &  \colorbox{green}{1} & \colorbox{light-blue}{2} \\
5 & 4 & 6 & \colorbox{yellow}{2} &  \colorbox{green}{3} & \colorbox{light-blue}{1} \\
\colorbox{yellow}{1} & \colorbox{yellow}{3} & \colorbox{yellow}{2} & 4 & 5 & 6 \\
\colorbox{green}{3} & \colorbox{green}{2} & \colorbox{green}{1}  &  5 & 6 & 4 \\
\colorbox{light-blue}{2} & \colorbox{light-blue}{1} & \colorbox{light-blue}{3} &  6 & 4 & 5 \\
\end{array}
\]

In more formalized manner construction which is described  in Examples \ref{Example_261} and  \ref{Example_262} is given in \cite{GONS_MARKOV_NECHAEV}. Using this construction  on base of T-quasigroups in \cite{GONS_MARKOV_NECHAEV} many MDS-codes are constructed. A pair of orthogonal quasigroups of order ten  is also constructed  \cite{GONS_MARKOV_NECHAEV, LIE_ZHU}.

\section{Belyavskaya modification}

G.B. Belyavskaya proposed modification of Bruck-Belousov construction \cite{B70_1, B70_4, B70_3}.

\begin{construction} \label{GBB_ALG} \index{prolongation!Belyavskaya}
We add additional column from the right  and additional row below, transpose  in these new cells all marked (transversal)  elements except one (in our example element ${\bf 2}$) in their fixed order and fill all remaining  empty cells except one  (with coordinates $(n+1, n+1)$) by the symbol "4". The cell with coordinates $(n+1, n+1)$ is filled by the not transposed  transversal element.
\[
{\displaystyle
\begin{array}{lllllll}
\begin{array}{lll}
1 & 2 &	\colorbox{yellow}{3}  \\
\colorbox{yellow}{2} &	3 &	1 \\
3 &	\colorbox{yellow}{1} &	2 \\
\end{array}
& \rightarrow &
\begin{array}{llll}
1 &	2 &	\Box &	\colorbox{yellow}{3} \\
\colorbox{yellow}{\bf 2}  &	3 &	1 & \Box \\
3& \Box	&	2 &	\colorbox{yellow}{1} \\
\Box &	\colorbox{yellow}{1} &	\colorbox{yellow}{3} &	\Box \\
\end{array}
& \rightarrow &
\begin{array}{llll}
1 &	2 &	4 &	\colorbox{yellow}{3} \\
\colorbox{yellow}{\bf 2} &	3 &	1 & 4 \\
3 &	4 &	2 &	\colorbox{yellow}{1} \\
4 &	\colorbox{yellow}{1} &	\colorbox{yellow}{3} &	\colorbox{yellow}{\bf 2} \\
\end{array}
\end{array}}
\]

%%
%% & \rightarrow &
%%\begin{array}{l|llll}
%%\star & 1 & 2 & 3 & 4\\
%%\hline
%%1&1 &	2 &	4 &	\colorbox{yellow}{3}\\
%%2&\colorbox{yellow}{\bf 2}  &	3 &	1 &	4\\
%%3&3&	4&	2&	\colorbox{yellow}{1}\\
%%4& 4 &	\colorbox{yellow}{1}&	\colorbox{yellow}{3} &	\colorbox{yellow}{\bf 2} \\
%%\end{array}
\end{construction}

\begin{construction} \label{CONSTRUCTION_1_GBB_0} Generalized Belyavskaya construction. It is possible to generalize Belyavskaya construction using more than one disjoin transversal.
\end{construction}

\begin{example}
We prolong the Latin square  of order $3$ (Example \ref{Latin_squere_Of_Ordr_3}) to  Latin square  of order $5$ using Belyavskaya prolongation construction (Construction \ref{GBB_ALG}) that simultaneously is applied to two transversals.
\end{example}

We add two additional columns from the right  and two additional rows below, transpose in these new cell all marked (transversal)  elements in their fixed order for exception of one element in any transversal. We take  element $1$ in yellow transversal and take element $1$ in green transversal. It is not obligatory that in  yellow and green transversal we take equal "exceptional" elements.

Step 1.

\[
\begin{array}{lllll}
\colorbox{yellow}{1} & \Box & \colorbox{light-blue}{3} & \Box &  \colorbox{green}{2}  \\
\colorbox{light-blue}{2} & \Box & \colorbox{green}{1} & \colorbox{yellow}{3} &  \Box   \\
\Box & \colorbox{light-blue}{1} & \Box & \colorbox{yellow}{2} &  \colorbox{green}{3}   \\
\Box & \colorbox{yellow}{3} & \colorbox{yellow}{2} &  \Box & \Box \\
\colorbox{green}{3} & \colorbox{green}{2} & \Box  &  \Box & \Box  \\
\end{array}
\]

Step 2.

Fill all remaining  empty after transposition transversal cells by the symbols "4, 5". In transversal cells of a fixed transversal we put the same element.

In the bottom of main diagonal write elements "4, 5".
Unfortunately here direct generalization of Belyavskaya construction is not possible.

\[
\begin{array}{lllll}
\colorbox{yellow}{1} & 4 & \colorbox{light-blue}{3} & 5 &  \colorbox{green}{2}  \\
\colorbox{light-blue}{2} & 5 & \colorbox{green}{1} & \colorbox{yellow}{3} &  4   \\
4 & \colorbox{light-blue}{1} & 5 & \colorbox{yellow}{2} &  \colorbox{green}{3}   \\
5 & \colorbox{yellow}{3} & \colorbox{yellow}{2} & 4 & 1 \\
\colorbox{green}{3} & \colorbox{green}{2} & 4  &  1 & 5  \\
\end{array}
\]

Modification of Step 1. We can change the order  colored  rows and/or  columns. For example we have changed 4-th and 5-th rows.

\[
\begin{array}{lllll}
\colorbox{yellow}{1} & 4 & \colorbox{light-blue}{3} & 5 &  \colorbox{green}{2}  \\
\colorbox{light-blue}{2} & 5 & \colorbox{green}{1} & \colorbox{yellow}{3} &  4   \\
4 & \colorbox{light-blue}{1} & 5 & \colorbox{yellow}{2} &  \colorbox{green}{3}   \\
\colorbox{green}{3} & \colorbox{green}{2} & 4  &  1 & 5  \\
5 & \colorbox{yellow}{3} & \colorbox{yellow}{2} & 4 & 1 \\
\end{array}
\]

\begin{example}
Using generalized Belyavskaya construction we  prolong the Latin square  of order $3$ to  Latin square of order $6$ in the following way.
\end{example}

Step 1.
We add three additional columns from the right  and three additional rows below, transpose in these new cells all marked (transversal)  elements in their fixed order for exception of the element $3$ from yellow transversal, the element $3$ from green transversal and the element $1$ from light-blue transversal.

\[
\begin{array}{llllll}
{\Box} & \Box & \Box & \colorbox{yellow}{1} &  \colorbox{green}{2} & \colorbox{light-blue}{3} \\
\Box & \colorbox{yellow}{3} & \Box & \Box &  \colorbox{green}{1} & \colorbox{light-blue}{2} \\
\colorbox{green}{3} & \colorbox{light-blue}{1} & {\Box} & \colorbox{yellow}{2} &  \Box & \Box  \\
\colorbox{yellow}{1} & \Box & \colorbox{yellow}{2} & \Box & \Box & \Box \\
\Box & \colorbox{green}{2} & \colorbox{green}{1}  &  \Box & \Box & \Box \\
\colorbox{light-blue}{2} &	\Box & \colorbox{light-blue}{3} &  \Box & \Box & \Box \\
\end{array}
\]

Step 2.
Fill all remaining  empty after transposition transversal cells by the symbols "4, 5, 6" in their  "old transversal order", i.e. we put the symbol $4$ in all empty cells of light-blue  transversal and so on.

\[
\begin{array}{llllll}
4 & 5 & 6 & \colorbox{yellow}{1} &  \colorbox{green}{2} & \colorbox{light-blue}{3} \\
6 & \colorbox{yellow}{3} & 5 & 4 &  \colorbox{green}{1} & \colorbox{light-blue}{2} \\
\colorbox{green}{3} & \colorbox{light-blue}{1} & 4 & \colorbox{yellow}{2} &  5 & 6  \\
\colorbox{yellow}{1} & 4 & \colorbox{yellow}{2} & \Box & \Box & \Box \\
5 & \colorbox{green}{2} & \colorbox{green}{1}  &  \Box & \Box & \Box \\
\colorbox{light-blue}{2} &	6 & \colorbox{light-blue}{3} &  \Box & \Box & \Box \\
\end{array}
\]

Step 3.
Remaining right bottom empty  square  we should complete in order to obtain a quasigroup. Bottom part of main diagonal we fill by the elements $3, 3, 1$, because namely these elements remain in transversals.

In this case cell $(5,4)$ can be filled only by the element 6. Remaining is clear for any fan of Sudoku.

\[
\begin{array}{llllll}
4 & 5 & 6 & \colorbox{yellow}{1} &  \colorbox{green}{2} & \colorbox{light-blue}{3} \\
6 & \colorbox{yellow}{3} & 5 & 4 &  \colorbox{green}{1} & \colorbox{light-blue}{2} \\
\colorbox{green}{3} & \colorbox{light-blue}{1} & 4 & \colorbox{yellow}{2} &  5 & 6  \\
\colorbox{yellow}{1} & 4 & \colorbox{yellow}{2} & 3 & 6 & 5 \\
5 & \colorbox{green}{2} & \colorbox{green}{1}  &  6 & 3 & 4 \\
\colorbox{light-blue}{2} &	6 & \colorbox{light-blue}{3} &  5 & 4 & 1 \\
\end{array}
\]

\section{Prolongation using quasicomplete mappings}

I.I.~Derienko and W.A.~Dudek propose prolongation construction of a quasigroup  using quasicomplete mappings \cite{DER_DUD}. This construction is generalization of Belyavskaya construction on quasicomplete mappings.

\begin{definition}
Let $(Q,\cdot)$ be a finite quasigroup,  $\sigma$ be a mapping of the set $Q$. We can construct the mapping $\overline{\sigma}$ in the following way:
\begin{equation}
\overline{\sigma} x = x \cdot \sigma x \quad \text{for all} \quad x\in Q.
\end{equation}
The mapping $\overline{\sigma}$ is called \textit{conjugated mapping}  to the mapping $\sigma$.
\end{definition}

A mapping $\sigma$ is \textit{quasicomplete}, if $\sigma$ is a permutation of a set $Q$ and   $\overline{\sigma}(Q)$
contains all elements of $Q$ except one. In this case there exists an element $a \in Q$,
called special, such that $a = \overline{\sigma} x_1  = \overline{\sigma} x_2 $
for some $x_1,  x_2 \in  Q$, $x_1 \neq x_2$.

\begin{construction} \index{prolongation!Derienko-Dudek}
 We start from a quasigroup and a quasicomplete mapping, act as in Belyavskaya construction but in the cell with coordinates $(n+1, n+1)$ we write the element $Q\setminus \overline{\sigma}Q$ (the special element).

\begin{example} \label{QUASI_COMPLET_EX}
We take the following quasigroup
\[
 \begin{array}{c|cccc}
\cdot & 1 & 2 & 3 & 4 \\
\hline
1 & \colorbox{yellow}{\bf{2}} & 1 & 3 & 4  \\
2 & 3 & 2 & \colorbox{yellow}{\bf{4}} & 1  \\
3 & 4 & \colorbox{yellow}{\bf{3}} & 1 & 2  \\
4 & 1 & 4 & 2 & \colorbox{yellow}{\bf{3}}
\end{array}
\]
and the following  mapping
\[
\sigma =
\left(
\begin{array} {llll}
1 & 2 & 3 & 4  \\
1 & 3 & 2 & 4   \\
\end{array}
\right).
\]
Then
\[
\overline{\sigma} =
\left(
\begin{array} {llll}
1 & 2 & 3 & 4  \\
2 & 4 & 3 & 3  \\
\end{array}
\right),
\]
$\sigma$ is quasicomplete mapping and $Q\setminus \overline{\sigma}Q= \{1\}$. Using Derienko-Dudek construction we obtain
\[
\begin{array}{llll}
\begin{array}{l|lllll}
\ast &  1 &  2 & 3 & 4 & 5  \\
\hline
1 & \Box & 1 & 3 & 4 & \bf{2} \\
2 & 3 & 2 & \Box & 1 & \bf{4} \\
3 & 4 & \Box  & 1 & 2 & \bf{3} \\
4 & 1 & 4 & 2 & \bf{3} &\Box \\
5 & \bf{2} & \bf{3} & \bf{4} & \Box &\Box \\
\end{array}
&\rightarrow &
\begin{array}{l|lllll}
\ast &  1 &  2 & 3 & 4 & 5  \\
\hline
1 & \colorbox{yellow}{5} & 1 & 3 & 4 & \bf{2} \\
2 & 3 & {2} & \colorbox{yellow}{5} & 1 & \bf{4} \\
3 & 4 & \colorbox{yellow}{5}  & {1} & 2 & \bf{3} \\
4 & 1 & 4 & 2 & {\bf{3}} & \colorbox{yellow}{5} \\
5 & \bf{2} & \bf{3} & \bf{4} & \colorbox{yellow}{5} & \colorbox{green}{1} \\
\end{array}
\end{array}
\]
\end{example}
\end{construction}

There exists a possibility to generalize Derienko-Dudek construction in Yamamoto spirit.
\begin{example}
It is clear that the following Latin square has four disjoint quasi-complete mappings. We shall use yellow and green quasi-complete mappings for the prolongation of this Latin square.
\[
 \begin{array}{cccc}
\colorbox{yellow}{2} & \colorbox{green}{1} & \colorbox{light-gray}{3} & \colorbox{light-blue}{4}  \\
\colorbox{green}{3} & \colorbox{light-blue}{2} & \colorbox{yellow}{4} & \colorbox{light-gray}{1}  \\
\colorbox{light-gray}{4} & \colorbox{yellow}{3} & \colorbox{light-blue}{1} & \colorbox{green}{2}  \\
\colorbox{light-blue}{1} & \colorbox{light-gray}{4} & \colorbox{green}{2} & \colorbox{yellow}{3}
\end{array}
\]

Step 1.

We add two columns and two rows and transpose there elements of the yellow and green quasi-complete mappings.

\[
\begin{array}{cccccc}
\Box & \Box  & \colorbox{light-gray}{3} & \colorbox{light-blue}{4} & \colorbox{yellow}{2}& \colorbox{green}{1} \\
\Box  & \colorbox{light-blue}{2} & \Box  & \colorbox{light-gray}{1} & \colorbox{yellow}{4}  &\colorbox{green}{3}\\
\colorbox{light-gray}{4} & \Box  & \colorbox{light-blue}{1} & \Box  & \colorbox{yellow}{3} & \colorbox{green}{2}  \\
\colorbox{light-blue}{1} & \colorbox{light-gray}{4} & \colorbox{green}{2} & \colorbox{yellow}{3} & \Box & \Box \\
\colorbox{yellow}{2} & \colorbox{yellow}{3} & \colorbox{yellow}{4} & \Box & \Box & \Box \\
\colorbox{green}{3} & \colorbox{green}{1} & \Box &\colorbox{green}{2} &  \Box & \Box \\
\end{array}
\]

Step 2.

We fill empty yellow cells in the initial  Latin square by the number $5$ and green cells by the number $6$.
We fill bottom part of main diagonal by the elements $1$ and $4$ respectively.

\[
\begin{array}{cccccc}
5 & 6  & \colorbox{light-gray}{3} & \colorbox{light-blue}{4} & \colorbox{yellow}{2}& \colorbox{green}{1} \\
6  & \colorbox{light-blue}{2} & 5  & \colorbox{light-gray}{1} & \colorbox{yellow}{4}  &\colorbox{green}{3}\\
\colorbox{light-gray}{4} & 5  & \colorbox{light-blue}{1} & 6  & \colorbox{yellow}{3} & \colorbox{green}{2}  \\
\colorbox{light-blue}{1} & \colorbox{light-gray}{4} & \colorbox{green}{2} & \colorbox{yellow}{3} & \Box & \Box \\
\colorbox{yellow}{2} & \colorbox{yellow}{3} & \colorbox{yellow}{4} & \Box & 1 & \Box \\
\colorbox{green}{3} & \colorbox{green}{1} & \Box &\colorbox{green}{2} &  \Box & 4 \\
\end{array}
\]

Step 3.

Finally we complement obtained partial Latin square to complete Latin square. It is easy to see that it is possible do this in a unique way.

\[
\begin{array}{cccccc}
5 & 6  & \colorbox{light-gray}{3} & \colorbox{light-blue}{4} & \colorbox{yellow}{2}& \colorbox{green}{1} \\
6  & \colorbox{light-blue}{2} & 5  & \colorbox{light-gray}{1} & \colorbox{yellow}{4}  &\colorbox{green}{3}\\
\colorbox{light-gray}{4} & 5  & \colorbox{light-blue}{1} & 6  & \colorbox{yellow}{3} & \colorbox{green}{2}  \\
\colorbox{light-blue}{1} & \colorbox{light-gray}{4} & \colorbox{green}{2} & \colorbox{yellow}{3} & 6 & 5 \\
\colorbox{yellow}{2} & \colorbox{yellow}{3} & \colorbox{yellow}{4} & 5 & 1 & 6 \\
\colorbox{green}{3} & \colorbox{green}{1} & 6 &\colorbox{green}{2} &  5 & 4 \\
\end{array}
\]

\end{example}

\section{Two step mixed procedure}

Suppose that a quasigroup of order $n$ has two disjoin transversals. On the first step   we can prolong this quasigroup to a qusigroup of order $n+1$ using Bruck-Belousov or Belyavskaya construction. After this procedure our  second transversal passes in complete or quasicomplete mapping (i.e. in $n$- or $(n-1)$-transversal) and we can prolong obtained quasigroup using  Bruck-Belousov or Belyavskaya construction, either Derienko-Dudek construction (if we have obtained quasicomplete mapping).

\begin{example}

\[
\begin{array}{lll}
\colorbox{yellow}{1} & \colorbox{green}{2} & \colorbox{light-blue}{3} \\
\colorbox{light-blue}{2} &\colorbox{yellow}{3} & \colorbox{green}{1} \\
\colorbox{green}{3} &	\colorbox{light-blue}{1} & \colorbox{yellow}{2} \\
\end{array}
\]
\end{example}

On the first step we use blue transversal and Belyavskaya construction (see Construction 2) and obtain the following Latin square.
\[
\begin{array}{llll}
1 &	2 &	\colorbox{light-blue}{4} &	{3} \\
\colorbox{light-blue}{\bf 2} &	3 &	1 & 4 \\
3 &	\colorbox{light-blue}{4} &	2 &	{1} \\
4 &	{1} &	{3} &	{\bf 2} \\
\end{array}
\]

In this case yellow transversal passes in $(n-1)$-thransversal (quasicomplete mapping)
\[
\sigma =
\left(
\begin{array} {llll}
1 & 2 & 3 & 4  \\
1 & 2 & 3 & 4   \\
\end{array}
\right).
\]
Then
\[
\overline{\sigma} =
\left(
\begin{array} {llll}
1 & 2 & 3 & 4  \\
1 & 3 & 2 & 2  \\
\end{array}
\right),
\]
$\sigma$ is quasicomplete mapping and $Q\setminus \overline{\sigma}Q= \{4\}$. Using Derienko-Dudek construction we obtain

\[
\begin{array}{lllll}
5 &	2 &	{4} &	{3} & 1 \\
{2} &	5 &	1 & 4 & 3 \\
3 &	{4} &	5 &	{1} & 2\\
4 &	{1} &	{3} &	{ 2} & 5 \\
1 &	{3} &	{2} &	{ 5} & 4\\
\end{array}
\]

\section{Contractions of quasigroups}

Procedure which is inverse to prolongation of quasigroups is called contraction of quasigroup. Procedures of contraction of quasigroups make from a quasigroup of order  $n$ a quasigroup of order $(n-1)$ or $(n-2)$. It is clear that any procedure of prolongation has its proper "inverse" procedure of contraction.  See \cite{B70_3, B70_4, DER_DER_DUD_13_II} for details.

Probably  procedures of prolongation and contraction of quasigroups can be used in cryptography.

\noindent \footnotesize{Institute of Mathematics and   Computer Science \\
Academy of Sciences of Moldova\\
 Academiei str. 5,  MD$-$2028 Chi\c{s}in\u{a}u \\
 Moldova\\
\noindent Email: scerb@math.md
}

\end{document}